\documentclass[12pt]{article}

\usepackage{amsmath}
\usepackage{amsthm}
\usepackage{amsfonts}
\usepackage{amssymb}

\newtheorem{thm}{Theorem}[section]
\newtheorem{defi}[thm]{Definition}
\newtheorem{prop}[thm]{Proposition}

\newtheorem{cor}[thm]{Corollary}

\newtheorem{prope}[thm]{Property}

\usepackage[dvips]{graphicx}

\begin{document}
\title{\bf Distance Hereditary Graphs and the Interlace Polynomial}
\author{Joanna A. Ellis-Monaghan$^1$\\{\small Department of Mathematics}\\
{\small Saint Michael's College}\\
{\small 1 Winooski Park, Colchester, VT 05439}\\
{\small Phone: (802)654-2660, Fax: (802)654-2610}\\
{\small jellis-monaghan@smcvt.edu}\\
\\Irasema Sarmiento\\
{\small Dipartimento di Matematica}\\
{\small  Universit\`a di Roma ``Tor Vergata"}\\
{\small Via della Ricerca Scientifica, I-00133, Rome, Italy}\\
{\small Phone: (39) 06 7259-4837, Fax: (39) 06 7259-4699}\\
{\small sarmient@mat.uniroma2.it}}\date{} \maketitle
\noindent
{\bf Proposed running head:} Distance Hereditary Graphs and the Interlace Polynomial.\\

$^1$ Research supported by the Vermont Genetics Network through the NIH Grant Number
1 P20 RR16462 from the BRIN/INBRE program of the NCRR.

\newpage

\begin{abstract}

\linespread{1.5}

The vertex-nullity interlace polynomial of a graph, described
by Arratia, Bollob\'as and Sorkin in
\cite{ABS00} as evolving from questions of DNA sequencing, and extended to
a two-variable interlace polynomial by the same authors
in \cite{ABSb}, evokes many open
questions. These include relations between the interlace polynomial
and the Tutte polynomial and the computational complexity of the
vertex-nullity interlace polynomial. Here, using the medial graph of
a planar graph, we relate the one-variable vertex-nullity interlace
polynomial to the classical Tutte polynomial when $x=y$, and
conclude that, like the Tutte polynomial, it is in general $\#$P-hard
to compute. We also show a relation between the two-variable interlace polynomial and the topological Tutte polynomial of
Bollob\'as and Riordan in \cite{BR01}.
We define the $\gamma $ invariant as the coefficient of
$x^1$ in the vertex-nullity interlace polynomial, analogously to the
$\beta$ invariant, which is the coefficient of $x^1$ in the Tutte
polynomial. We then turn to distance hereditary graphs, characterized
by Bandelt and Mulder in
\cite{BM86} as being constructed by a sequence of adding pendant and twin
vertices, and show that graphs in this class have $\gamma$ invariant of
$2^{n+1}$ when $n$ true twins are added in their construction. We furthermore
show that bipartite distance hereditary graphs are exactly the class
of graphs with $\gamma$ invariant $2$, just as the series-parallel graphs
are exactly the class of graphs with $\beta$ invariant $1$. In addition,
we show that a bipartite distance hereditary graph arises precisely
as the circle graph of any Euler circuit in the oriented medial
graph of a series-parallel graph. From this we conclude that the
vertex-nullity interlace polynomial is polynomial time to compute for
bipartite distance hereditry graphs, just as the Tutte
polynomial is polynomial time to compute for series-parallel graphs.\\



\end{abstract}
\footnotetext{DH Graphs and the Interlace Polynomial.} \newpage

\linespread{1}

\section{Introduction}\label{introduction}

In \cite{ABS00}, Arratia, Bollob\'as and Sorkin
defined a one-variable graph polynomial
$q_N(G)$ (denoted $q$ there, but we follow their later work \cite{ABSb},
reserving $q$ for the two-variable generalization) motivated by questions
arising from DNA sequencing by hybridization
addressed by Arratia, Bollob\'as, Coppersmith and Sorkin in
\cite{ABCS00}. This
polynomial models the interlaced repeated subsequences of DNA that can
interfere with the unique reconstruction of the original DNA strand.
This work promptly generated further interest and other applications
by Arratia, Bollob\'as, and Sorkin  \cite{ABS}, \cite{ABSb},
Aigner and van der Holst \cite{AH04},
Ballister, Bollob\'as, Cutler, and Pebody \cite{BBCP02},
and Ballister, Bollob\'as, Riordan, and Scott \cite{BBRS01}. In \cite{ABSb},
Arratia, Bollob\'as, and Sorkin define a
two-variable interlace polynomial, and show that the original polynomial of
\cite{ABS00} is a specialization of it, renaming the original interlace
polynomial the vertex-nullity interlace polynomial due to its
relationship with the two-variable generalization. In \cite{ABS00}, and
again in \cite{ABS}, Arratia, Bollob\'as, and Sorkin provide some tantalizing properties of
the vertex-nullity interlace polynomial, particularly for
circle graphs, but describe it as ``...territory newly invented,
but not yet explored." They present several open questions, some
reiterated in \cite{ABSb} in  the context of the two-variable interlace
polynomial, including the interlace polynomial's relation to known graph
polynomials and the computational complexity of the vertex-nullity interlace
polynomial.

\footnotetext{DH Graphs and the Interlace Polynomial.}

In Section \ref{toptutte} we note that, for planar graphs, the vertex-nullity interlace polynomial
is related to the Tutte polynomial via medial graphs, and from this
it follows that the vertex-nullity interlace polynomial, like the
Tutte polynomial, is computationally intractable in the general case.
We similarly provide a relationship between the
two-variable interlace polynomial and the topological Tutte
polynomial of Bollob\'as and Riordan from \cite{BR01}, again via
a medial graph construction, thus showing that the interlace
polynomial in some sense is encoding topological information.

After establishing that a polynomial is $\#$P-hard to compute in
general, a natural question is whether there are interesting classes of
graphs for which it is tractable. We identify such a class of graphs,
bipartite distance hereditary graphs,
for which the vertex-nullity interlace polynomial is
polynomial time to compute. In Section \ref{last} we show that these graphs are
characterized by an invariant $\gamma$ (the coefficient of $x^1$
in $q_N$) analogously to the way series-parallel graphs are characterized
by the $\beta$ invariant (the common coefficient of $x^1$ and $y^1$
in the Tutte polynomial). This completes the characterization of graphs
for which the coefficient of $x^1$ in $q_N$ is $2$, initiated
by Aigner and van der Holst in
\cite{AH04}.

We show that the $\gamma$ invariant of a (not necessarily bipartite)
distance hereditary graph $G$ is $2^{n+1}$ if $n$ true twins are
added in some construction sequence of $G$, but that distance
hereditary graphs do not comprise the entire class for which
$\gamma$ is a power of $2$.

We  show in Section \ref{PD} that bipartite distance hereditary
graphs arise precisely as circle graphs derived from
Eulerian circuits in the oriented medial graphs of series-parallel
graphs. From this characterization, and that distance
hereditary graphs may be recognized in polynomial time, we conclude
that the vertex-nullity interlace polynomial is polynomial time to compute
on the class of bipartite distance hereditary graphs.

\section{The interlace and circuit partition polynomials}\label{intercircuit}
\label{complexity}

The vertex-nullity interlace polynomial of a graph was defined
recursively by Arratia, Bollob\'as, and Sorkin in \cite{ABS00}
via a pivoting operation and was seen by them in \cite{ABSb} to be
a specialization of a much richer two-variable interlace
polynomial, $q(G;x,y)$, with a similar
pivot recursion.  Let $vw$ be an edge of a graph $G$, and let $A_v$,
$A_w$ and $A_{vw}$ be the sets of vertices of $G$ adjacent to $v$ only, $w$ only,
and to both $v$ and $w$, respectively.  The pivot operation ``toggles"
the edges between $A_v$,  $A_w$ and $A_{vw}$, by deleting existing edges and inserting
edges between previously non-adjacent vertices.
The result of this operation
is denoted $G^{vw}$.  More formally,  $G^{vw}$  has the same vertex set as $G$, and edge
set equal to the symmetric difference $E(G)\Delta   S$, where $S$  is the
complete tripartite graph with vertex classes
$A_v$,  $A_w$ and $A_{vw}$. See Figure $1$.\\

\footnotetext{DH Graphs and the Interlace Polynomial.}

[Insert Figure 1: The pivot operation.]

\begin{defi}\label{vertexnullity}

The vertex-nullity interlace polynomial is defined recursively as:

$$
q_N(G;x)=
\begin{cases}
x^n {\text{ if }}G=E_n, {\text{ the edgeless graph on $n$ vertices}}\\
q_N(G- v;x)+q_N(G^{vw}- w;x) {\text{ if }}vw\in E(G)
\end{cases}
$$

\end{defi}

This polynomial was shown to be well-defined on all simple graphs in \cite{ABS00}.

\begin{defi}

The two-variable interlace polynomial is defined,
for a graph $G$ of order $n$,
by

\begin{equation}\label{redq}q(G;x,y)=\sum_{S\subseteq V(G)}(x-1)^{r(G[S])}(y-1)^{n(G[S])},\end{equation}

\noindent
where $r(G[S])$ and $n(G[S])=|S|-r(G[S])$ are, respectively,
the ${\mathbb F}_2$-rank and
nullity of the adjacency
matrix of $G[S]$, the subgraph of $G$ induced by $S$.

\end{defi}

Equivalently, the two-variable interlace polynomial can be defined by the
following reduction formulas from Arratia, Bollob\'as and Sorkin
\cite{ABSb}.

For a graph $G$, for any edge $ab$ where neither $a$ nor $b$ has
a loop,

\begin{equation}\label{2eq}
q(G)=q(G- a)+q(G^{ab}- b)+((x-1)^2-1)q(G^{ab}- a- b), \\
\end{equation}

for any looped vertex $a$,

$$q(G)=q(G- a)+(x-1)q(G^a- a),$$

\noindent
and, for the edgeless graph $E_n$ on $n\geq 0$ vertices, $q(E_n)=y^n$.
Here $G^a$ is the {\em local complementation} of $G$ , and is defined
as follows. Let $N(a)$ be the neighbors of $a$, that is, the set
$\{w\in V: a {\text{ and }}w {\text{ are joined by an edge}}\}$.
Thus $a\in N (a)$ iff $a$ is a loop.
The graph $G^a$ is equal to $G$ except
that $G^a[N (a)] =\overline{G[N (a)]}$, i.e.
we ``toggle" the edges among the neighbors of $a$,
switching edges to non-edges and vice-versa.

Arratia, Bollob\'as, and Sorkin show in \cite{ABSb}
that the vertex-nullity interlace polynomial is a
specialization of the two-variable interlace
polynomial as follows:

\footnotetext{DH Graphs and the Interlace Polynomial.}

\begin{equation}\label{stateqN}
q_N(G;y)=q(G;2,y)=\sum_{W\subseteq V(G)}(y-1)^{n(G[W])}.
\end{equation}

An equivalent formulation for $q_N(G;x)$  is given by Aigner and van der Holst in
\cite{AH04}.


\begin{prop}\label{5.8'}
For any graph $G$ with more than one vertex, the coefficients of $x^1$ and $y^1$ in the
two-variable
interlace polynomial are one the negative of the other.
\end{prop}
\begin{proof}
We collect terms and write $q(G;x,y)=\sum_{i,j\geq 0}a_{ij}x^iy^j$ for the
two-variable interlace
polynomial.
Thus, we want to prove that $a_{10}=-a_{01}$. That
is, that $\frac{\partial{q}}{\partial x}|_{x=y=0}=a_{10}=-a_{01}=-
\frac{\partial{q}}{\partial y}|_{x=y=0}$.
Note that
$\frac{\partial{q}}{\partial x}|_{x=y=0}$

$=\sum_{S\subset V(G)}
(-1)^{n(G[S])+r(G[S])-1}r(G[S])$,

and $\frac{\partial }{\partial y}
|_{x=y=0}=\sum_{S\subseteq V(G)}(-1)^{n(G[S])+r(G[S])-1}n(G[S])$.

Thus
$a_{10}+a_{01}=\sum_{S\subset V(G)}
(-1)^{n(G[S])+r(G[S])-1}(r(G[S])+n(G[S]))$

$=-\sum_{S\subseteq V(G)}(-1)^{|S|}|S|$
$=-\sum_{l=0}^{|V(G)|}(-1)^l\binom{|V(G)|}{l}l$
$=(-1)^{|V(G)|+1}\frac{d}{dt}(t-1)^{|V(G)|}|_{t=1}=0$.

Therefore $a_{10}+a_{01}=0$, so
$a_{10}=-a_{01}.$

\end{proof}

\footnotetext{DH Graphs and the Interlace Polynomial.}

Although we do not pursue it here, we would not be surprised by
relations among the coefficients of $q$ analogous to those for the
Tutte polynomial found by Brylawski  in \cite{Bry80}.

\begin{cor}\label{2.4}

 If $q_N(G;x)=\sum a_ix^i$, then $a_1=\sum_i a_{i,1}2^i$, so
 $a_1=a_{01}=-a_{10}$ if and only if $\sum_{i\geq 1} a_{i,1}2^i=0$.

\end{cor}
\begin{proof}

 $\sum_ia_ix^i=q_N(G;x)=q(G;2,x)=\sum_{i,j}a_{i,j} 2^ix^j$, which
 with Proposition \ref{5.8'}   gives the result.

\end{proof}

\begin{prop}\label{pos}

Let $G$ be a simple graph. Then the non-zero coefficients of
$q_N(G;x)$ are positive integers.

\end{prop}

\begin{proof}

The proof is by induction on the number of vertices of $G$.
If $G$ has one vertex and no edges, then
$q_N(G;x)=x$.
Now assume the hypothesis holds for all
simple graphs on $n-1$ vertices, and let $G$ be a simple
graph with $n$ vertices. If $G$ has no edges,
$q_N(G)=x^n$. Otherwise, let $uv$ be an edge of $G$. Thus,
$q_N(G;x)=q_N(G- u;x)+q_n(G^{uv}- v;x)$.
By induction, the
non-zero coefficients of
$q_N(G- u;x)$ and
$q_N(G^{uv}- v;x)$ are
positive integers, and thus the non-zero
coefficients of
$q_N(G;x)$
are sums of non-negative integers.

\end{proof}

In \cite{ABS00} and \cite{ABSb}, Arratia, Bollob\'as, and
Sorkin  give an interpretation
of the vertex-nullity interlace polynomial of a
circle graph in terms of the circuit partition, or Martin, polynomial of a related
$4$-regular Eulerian digraph.  Recall that a {\em circle graph}
on $n$ vertices is a graph $G$
derived from a chord diagram,
where two copies of each of the symbols $1$ through $n$ are
arranged on the perimeter of a circle,
and a chord is drawn between like symbols.
Two vertices
$v$ and $w$ in $G$ share an edge if and only if their corresponding
chords intersect in the chord diagram.
See Figure $2$.

Circle graphs have also been called alternance graphs by
Bouchet \cite{Bou88} and
interlace graphs by Arratia, Bollob\'as and Sorkin \cite{ABS00}.
Research on circle graphs includes a
complete characterization and a polynomial time algorithm for
identifying them. For example see
Bouchet \cite{Bou85}, \cite{Bou87b}, \cite{Bou87c}, \cite{Bou94},
Czemerinski, Dur\'an, and Gravano \cite{CDG02},
Dur\'an \cite{Dur03}, Fraysseix \cite{Fra84},
Gasse \cite{Gas97},
Read and  Rosenstiehl \cite{RR78a}, \cite{RR78b},
and Wessel and P\"oschel \cite{WP84}.\\

[Insert Figure 2: The circle graph of a chord diagram.]\\

\footnotetext{DH Graphs and the Interlace Polynomial.}

A $4$-{\em regular Eulerian} digraph is a $4$-regular directed graph such that,
at each vertex, two edges are oriented inward, and two are oriented outward.
A $4$-regular Eulerian digraph is called a $2$-in, $2$-out graph in \cite{ABS00}. Note that if
$C$ is an Eulerian circuit of a $4$-regular Eulerian digraph, and we write the vertices
along the perimeter of a circle in the order that they are visited by $C$ (each is visited
exactly $2$ times), and then draw a chord between like vertices, the result is a chord diagram.

\begin{defi}
\rm
A {\em graph state} of a $4$-regular Eulerian digraph $\vec{G}$  is the result of replacing
each $4$-valent vertex $v$ of  $\vec{G}$ with two $2$-valent vertices each
joining an incoming and an outgoing
edge originally adjacent to $v$.  Thus a graph state is a disjoint union of consistently
oriented cycles.  See Figure $3$.
\end{defi}

Note that graph states (see \cite{E-M98}) are equivalent to the circuit
partitions of Arratia, Bollob\'as and Sorkin \cite{ABS00}
and Bollob\'as \cite{Bol02}, the Eulerian
decompositions of Bouchet \cite{Bou88}, and the Eulerian
$k$-partitions of Martin \cite{Mar77} and Las Vergnas \cite{Las83}.\\

[Insert Figure $3$: A graph state.]\\

\begin{defi}
\rm
The {\em circuit partition polynomial} of a $4$-regular Eulerian
digraph  $\bar{G}$ is
$f(\vec{G};x)=\sum_{k\geq 0}f_k(\vec{G})x^k$, where $f_k(\vec{G})$
is the number of graph states of $\bar{G}$  with $k$ components, defining
$f_0(\vec{G})$ to be $1$ if $\vec{G}$  has no edges,
and $0$ otherwise.
\end{defi}

The circuit partition polynomial is a simple translation
of the Martin polynomial $m(\vec{G};x)$,
defined
recursively for $4$-regular digraphs by Martin in his 1977 thesis
\cite{Mar77}, with
$f(\vec{G};x)=xm(\vec{G};x+1)$.

Las Vergnas found closed forms for the Martin polynomials (for both
graphs and digraphs). He also extended their properties to general
Eulerian digraphs and further developed their theory
(see \cite{Las79}, \cite{Las88}, \cite{Las83}). The transforms of the
Martin polynomials, for arbitrary Eulerian graphs and digraphs,
were given in \cite{E-M98}, and then aptly named circuit partition
polynomials by Bollob\'as in \cite{Bol02}, with splitting
identities provided in \cite{Bol02} and \cite{E-Mb}. The circuit partition
polynomial is also a specialization of a much broader
multivariable polynomial, the generalized transition polynomial of
\cite{E-MS02}, which assimilates such graph invariants as the
Penrose polynomial that are not evaluations of the
Tutte polynomial.

For circle graphs, the vertex-nullity
interlace polynomial and the circuit partition polynomial
are related
by the following theorem.

\begin{thm}(Arratia, Bollob\'as and Sorkin \cite{ABS00}, Theorem 6.1).\label{A}

If $\vec{G}$  is a $4$-regular Eulerian digraph, $C$ is any Eulerian
circuit of $\vec{G}$ , and  $H$ is the circle graph of the chord diagram
determined by $C$,
 then $f(\vec{G};x)=xq_N(H;x+1)$ .

\end{thm}

\footnotetext{DH Graphs and the Interlace Polynomial.}

\section{Relation to the classical and topological
Tutte polynomials and computational complexity}\label{toptutte}

Theorem \ref{A}, combined with a relationship between the Martin and Tutte polynomials for
planar graphs, relates the vertex-nullity interlace polynomial to the Tutte polynomial and consequently
resolves the computational complexity question
raised by Arratia, Bollob\'as and Sorkin in \cite{ABS00} and  \cite{ABSb}.
Arratia, Bollob\'as and Sorkin prove in \cite{ABSb} that the two-variable interlace and the
one-variable vertex-rank polynomials are $\# P$-hard to compute,
with only  the computational
complexity of the original vertex-nullity polynomial  left unresolved.

Let $G$ be a plane graph.
Its medial graph, $G_m$, has vertices corresponding to
the edges of $G$.
Two vertices of $G_m$ are joined by an edge if
the corresponding edges of $G$ are neighbors in the
cyclic order around a vertex.
We then color the faces of the medial graph black
or white, depending on whether they contain or do not contain, respectively, a vertex of the
original graph $G$.  This face-$2$-colors the medial graph.  The edges of the medial graph are
then directed so that the black face is on the left of an incident
edge.  See Figure $4$. Denote this oriented
medial graph by $\vec{G}_m$.\\

[Insert Figure $4$: The medial graph.]\\

The {\em Tutte polynomial} of a graph, $t(G;x,y)$,
may be defined by the linear recursion relation
$t(G;x,y)=t(G- e;x,y)+t(G/e;x,y)$
if $G$ has an edge $e$ that is neither an isthmus
(cut-edge or bridge) nor a loop of $G$, and by
$t(G;x,y)=x^iy^j$ if $G$ consists of $i$
isthmuses and $j$ loops. See Brylawski \cite{Bry80} or
Brylawski and Oxley \cite{BO92}, for example, for an
in-depth treatment of the
Tutte polynomial, including generalizations to matroids.

Martin (\cite{Mar77}, \cite{Mar78}) found the relationship
$m(\vec{G}_m;x)=t(G;x,x)$, which was further explored by
Las Vergnas in \cite{Las79}, \cite{Las88}.  This now allows
us to relate the vertex-nullity interlace polynomial
to the Tutte polynomial,
a relation also observed by Arratia, Bollob\'as and Sorkin at the end of
Section $7$ in \cite{ABS}.

\begin{thm}\label{B}

If $G$ is a planar graph, and $H$ is the circle graph of some Eulerian
circuit of  $\vec{G}_m$,
then $q_N(H;x)=t(G;x,x)$.

\end{thm}

\begin{proof}
By Theorem \ref{A}, $f(\vec{G}_m;x)=xq_N(H;x+1)$, but
recalling that the circuit partition and Martin polynomials are
simple translations of each other,
we have that $f(\vec{G}_m;x)=xm(\vec{G}_m;x+1)$, and hence
$q_N(H;x)=m(\vec{G}_m;x)=t(G;x,x)$.
\end{proof}

\footnotetext{DH Graphs and the Interlace Polynomial.}

Since the Tutte polynomial is known to be $\#$P-hard for planar graphs except at the
isolated points $(1,1)$, $(-1, -1)$, $(j, j^2)$, $(j^2,j)$
(where  $j=e^{\frac{2\pi i}{3}}$),
and along the curves $(x-1)(y-1)=1$ and
$(x-1)(y-1)=2$,
(see Jaeger, Vertigan, and Welsh \cite{JVW90} and
Welsh \cite{Wel93}), we have the following immediate corollary.

\begin{cor}\label{Ca}
The vertex-nullity interlace polynomial is $\#$P-hard in general.
\end{cor}

Of course, this leads immediately to the question of whether there might be
classes of graphs for which the vertex-nullity
interlace polynomial might be more tractable.
We provide one such class in Sections \ref{last} and \ref{PD}.

A natural question arises as to whether
Theorem \ref{B} might be extended. The idea may indeed be applied
elsewhere, and interestingly,
to the topological Tutte polynomial of Bollob\'as and
Riordan (\cite{BR01}), where the classical
Tutte polynomial is generalized to encode topological information about
graphs embedded on orientable surfaces.

The topological Tutte polynomial of Bollob\'as and Riordan for cyclic graphs
was defined in \cite{BR01}, and generalized for non-orientable
ribbon graphs by Bollob\'as and Riordan in \cite{BR02}.
{\em Cyclic graphs}
are graphs with {\em rotation systems}, that
is a family of local rotations around each vertex of
$G$. A {\em local rotation} around a vertex $v$ is a cyclic
order of the edges incident with $v$.
Contraction and deletion for cyclic graphs are
described by Bollob\'as and Riordan in \cite{BR01}. If $e$ is an edge of
a cyclic graph
${\bf G}$, then ${\bf G}- e$ is the cyclic graph obtained
by  deleting the edge $e$ from the underlying graph and from
whichever local rotations in which it occurs.
The contraction ${\bf G}/e$ of a non loop edge
$e=uv$ has $G/e$ as its underlying graph. Let $w$ be the vertex
of $G/e$ obtained by identifying $u$ and $v$.
The local rotation at $w$ is obtained by uniting those at
$u$ and $v$ using $e$. That is, following the edges after $e$ in the
local rotation at $u$ until we get to $e$ again. From there,
we follow the edges after $e$ in the local rotation around $v$.

A cyclic graph ${\bf G}$ with a single
vertex $v$ is given by the
cyclic order of the half edges around $v$. Therefore we can
identify ${\bf G}$ with the chord diagram $D$ that has
labels corresponding to the edges of ${\bf G}$ around
the boundary in exactly the same
order given by the cyclic permutation of the
edges around $v$. Cyclic graphs correspond to graphs embedded in
oriented surfaces (see Bollob\'as and Riordan
\cite{BR01} and the references therein).

\footnotetext{DH  Graphs and the Interlace Polynomial.}

 The ${\mathbb F}_2$-rank of a chord diagram $D$ was defined
 by Bollob\'as and Riordan in \cite{BR01} as
 $r(D)=\frac{1}{2} r(M(H))$, where $M(H)$ is the adjacency matrix
 of the circle graph $H$ of $D$, and
 $n(D)$ is the number of chords. It is observed that if
 ${\bf G}$ is a cyclic graph with one vertex, then
 $C({\bf G};X,Y,Z)=\sum_{D'\subseteq D}Y^{n(D')}Z^{r(D')}$.
 Here $D$ is the chord diagram determined by the cyclic order of
 the loop half-edges about the single vertex of ${\bf G}$, and
 $D'$ is a subchord diagram (the chord diagram formed from a
 subset of the chords of $D$). ${\bf G}$ and $D$ determine each other
 up to isomorphism.

\begin{defi}

Let ${\bf G}$ be a cyclic graph. The polynomial
$C({\bf G};X,Y,Z)$ is defined by

$$C({\bf G};X,Y,Z)=C({\bf G}/e;X,Y,Z)+C({\bf G}- e;X,Y,Z)$$

\noindent
if
$e$ is neither a bridge nor a loop of ${\bf G}$, and
$C({\bf G};X,Y,Z)=XC({\bf G}/e;X,Y,Z)$ if $e$ is a bridge in
${\bf G}$.

For one vertex cyclic graphs we have
$C(D;X,Y,Z)=\sum_{D'\subseteq D}Y^{n{D'}}Z^{r(D')}$, where $D$ is the
chord diagram corresponding to ${\bf G}$.

\end{defi}

\begin{thm}\label{Cpoly}

 Let $H$ be the circle graph of a chord diagram $D$. Then

 $$q(H;YZ^{\frac{1}{2}}+1,Y+1)=C(D;X,Y,Z),$$

 and

 $$q(H;x,y)=C(D;x,y-1,(\frac{x-1}{y-1})^2).$$

\end{thm}
\begin{proof}

 The vertices of $H$ correspond to the chords of $D$. A set
 $S\subseteq V(H)$ corresponds to a subdiagram $D_S$ of $D$.
 The induced subgraph $H[S]$ is the circle graph of
 $D_S$.

\footnotetext{DH Graphs and the Interlace Polynomial.}

 It follows that $r(D_S)=\frac{1}{2}r(H[S])$. Thus

 $$C(D;X,Y,Z)=\sum_{D'\subseteq D}Y^{n(D')}Z^{r(D')}$$

 $$=\sum_{S\subseteq V(H)}Y^{|S|}Z^{\frac{1}{2}r(H[S])}$$

 $$=\sum_{S\subseteq V(H)} (YZ^{\frac{1}{2}})^{r(H[S])}Y^{n(H[S])}$$

 $$=q(H;YZ^{\frac{1}{2}}+1,Y+1).$$

 Therefore,

 $$ q(H;YZ^{\frac{1}{2}}+1,Y+1)=C(D;X,Y,Z), {\text{ and }}$$

 $$q(H;x,y)=C(D;x,y-1,(\frac{x-1}{y-1})^2).$$

\end{proof}

 By Theorem  \ref{Cpoly} we obtain a reduction
 formula for the $C$-polynomial of a chord diagram.

\footnotetext{DH Graphs and the Interlace Polynomial.}

\begin{defi}

Let $D$ be a chord diagram and $a$, $b$ be two intersecting
chords. Let $H$ be the circle graph of $D$. Then $D^{ab}$ is the chord
diagram whose circle graph is $H^{ab}$.

\end{defi}

\begin{cor}

 Let $D$ be a chord diagram and let $a$, $b$ be two intersecting chords.
Then

\noindent
$C(D;X, Y,Z)=C(D- a;X, Y,Z)+C(D^{ab}- a;X, Y,Z)
+(Y^2Z-1)$
\noindent
$C(D^{ab}- a- b;X, Y,Z).$

\end{cor}

The medial graph, for graphs embedded in surfaces, is defined
as in the plane case.
Loops homotope to zero or bridges of ${\bf G}$ are cut vertices of $\vec{{\bf G}}_m$.
Thus they correspond to isolated
chords in the
chord diagram $D$ of any Eulerian circuit of $\vec{{\bf G}}_m$ and to isolated
chords in $H$. If $H$ is
the circle graph of $D$, then by Theorem \ref{Cpoly},
$q(H;x,y)=C(D;x,y-1,(\frac{x-1}{y-1})^2)$.
But if $e$ is a loop homotope to zero or a bridge,
$D=D_1+D_2$ where $D_1$ is the isolated chord corresponding to
$e$. By Corollary $2$ in Bollob\'as and Riordan \cite{BR01}
we have that $C({\bf {\bf G}})=C(D_1)C(D_2)$.
Thus $q(H;x,y)=C(D_1;x,y-1,(\frac{x-1}{y-1})^2)
C(D_2;x,y-1,(\frac{x-1}{y-1})^2)=y q(H[V- e];x,y)$
and this is independent of whether $e$ is a loop homotope to zero or a bridge.
Therefore, if $e$ is a loop homotope to zero or a bridge,
$q(H;x,y)=yq(H[V- e];x,y)$. It follows that the
$q(H;x,y)$ does not distinguish loops homotope to zero and bridges of $G$. Thus, in the
general case, $q(H;x,y)$ is not an evaluation of the
Tutte polynomial of $G$.

\section{Distance hereditary graphs and the $\gamma$ invariant}\label{last}

We turn our attention now to the question raised in Section
\ref{toptutte}
about classes of graphs for which the vertex-nullity
interlace polynomial is
polynomial time computable. Series-parallel graphs
and their characterization
by the $\beta$ invariant lead to the considerations of this
section. Recall that a series-parallel graph is (necessarily)
a planar graph
constructed from a digon  by repeatedly adding an edge in
parallel to an existing edge by including a multiple edge, or
adding an edge in series with an existing edge
by inserting a vertex of degree $2$ into the edge.
Also recall that for a graph with $2$ or more edges, the coefficients of
$x^1$ and $y^1$ in the Tutte polynomial are equal, and this common value,
$\beta (G)$, introduced by Crapo in \cite{Cra67}, is called
the $\beta$ invariant of $G$.
Brylawski characterized series-parallel graphs
(in the more general context of matroids) in
\cite{Bry71} by the property that $G$ is a
series-parallel graph if and only if $\beta (G)=1$.
The $\beta$ invariant has been explored further, for example by
Oxley in
\cite{Oxl82} and by
Benashki, Martin, Moore, and  Traldi  in   \cite{BMMT95}.
In analogy with the $\beta$ invariant, we define the $\gamma$ {\em invariant}
as the coefficient of $x^1$ in $q_N(G;x)$. This section
examines distance hereditary graphs ({\em DH graphs}) and characterizes
them in terms of the $\gamma$ invariant in that
$\gamma (G)=2^n$ if $G$ is a DH graph,
and $\gamma=2$ if and only if $G$ is a bipartite DH
graph ({\em BDH graph}).
En route, we present several general reduction formulas
for $\gamma (G)$, $q_N(G;x)$ and $q(G;x,y)$.

We show in Proposition \ref{5.20} that $H$ is a BDH graph if and only if it is
a circle graph corresponding to an Eulerian cycle in the
medial graph of some series-parallel graph,
and this
suffices
to show $\gamma (G)=2$ if $G$ is a BDH graph, but not vice versa. In
Section \ref{PD} we show that $q_N$ is polynomial time to compute for
BDH graphs.

\begin{defi}
\rm
If $q_N(G;x)=\sum a_i x^i$, then we call $a_1$ the $\gamma$ {\em invariant}
of $G$
and denote it $\gamma (G)$.
\end{defi}

We recall from Proposition 5.3 in
Arratia, Bollob\'as, and Sorkin  \cite{ABS00} that $q_N$ is multiplicative
on disjoint unions, and then note some very simple properties of
$\gamma$ that will be used repeatedly throughout this section.

\begin{prope}\label{5.2}
If $vw$ is an edge of $G$, then
$\gamma (G)=\gamma (G- v)+\gamma (G^{vw} - w)$.
\end{prope}

\footnotetext{DH Graphs and the Interlace Polynomial.}

\begin{prope}\label{5.2a}
If $vw$ is an edge of $G$, then $\gamma (G)=\gamma (G^{vw})$.
\end{prope}
\begin{proof}

The property is immediate from  $q_N(G;x)=q_N(G^{uv};x)$,
proposition $8.1$ in
Arratia, Bollob\'as, and Sorkin   \cite{ABS00}.

\end{proof}

\begin{prope}\label{5.3}
$\gamma (G)=0$ if and only if $G$ has more than
one component.
\end{prope}
\begin{proof}

The property is immediate from
Remark $20$ of
Arratia, Bollob\'as, and  Sorkin   \cite{ABS} which states that the
degree of the lowest degree term of $q_N(G)$ is the number of components
of $G$.
\end{proof}

\begin{prope}\label{5.4}
If $G$ is a connected graph, then $\gamma (G- v)=0$ if and only
if $v$ is a cut vertex of $G$ or $G$ consists of a single vertex.
\end{prope}

\begin{prope}\label{5.5}
If $G$ is a connected graph, then $\gamma (G)=1$ if and
only if $G$ is an isolated vertex.
\end{prope}
\begin{proof}
That $\gamma (G)\neq 1$ if $|V(G)|\geq 2$ follows from Property
\ref{5.2} and induction on $|E(G)|$, since both
$\gamma (G- v)$ and $\gamma (G^{uv}- u)$ are either $0$ or
greater than $1$, unless $G=K_2$, in which case both are $1$.
If $G$ is an isolated vertex, then $\gamma (G)=1$.
\end{proof}

We also note that we immediately derive a new interpretation for the
$\beta$ invariant of a planar graph.

\begin{prop}\label{5.6}
Let $G$ be a planar graph and let $H$ be the circle graph of some
Eulerian circuit of $\vec{G}_m$. Then

$$\beta (G)=-\frac{1}{2}\sum_{S\subseteq V(G)}n(G[S])(-1)^{n(G[S])}
=\frac{1}{2}a_1=\sum_i a_{i,1}2^{i-1}.$$

\end{prop}
\begin{proof}
By Theorem \ref{B}, $q_N(H;x)=t(G;x,x)$, so
$\gamma (H)=2\beta (G)$. But
$\gamma (H)=q_N'(H;0)$, so the result follows from
the expression for $q_N$ in equation (\ref{stateqN}) and Corollary \ref{2.4}.

\end{proof}

\footnotetext{DH Graphs and the Interlace Polynomial.}

We now review distance hereditary graphs. This important class of graphs,
introduced by E. Howorka, has a number of characterizations
(see \cite{How77a}, \cite{How77b}), the one
responsible for the name being that any connected induced subgraph of a
distance hereditary graph $G$ inherits its distance function from
$G$. See Brandst\"adt, Le, and Spinrad \cite{BLS99} and
McKee and McMorris \cite{MM99} for detailed overviews of this
and related classes of graphs.

However, for our purposes, we will use the formulation of
Definition \ref{5.7}, due to Bandelt and Mulder in \cite{BM86},
although we do not allow
infinite graphs. Twin vertices (see Burlet and Uhry's usage in
\cite{BU82}) are
``split pairs" in \cite{BM86}; we use the terminology twin because
of the emphasis on the differing roles of strong vs. weak twins with
respect to the interlace polynomial. Also, restricting
twins to non-isolated vertices and starting with $K_1$, as we
have in Definition \ref{5.7}, is equivalent to starting with
$K_2$ as in \cite{BM86}, both formulations serving to assure that the
resulting graph is connected. In fact, our restriction forces $K_2$ to
always be the next step after $K_1$.

\footnotetext{DH Graphs and the Interlace Polynomial.}

\begin{defi}\label{5.7}
\rm
A {\em distance hereditary } graph (DH graph) is a graph
that can be constructed from a single vertex by a finite number of
applications of
the following operations:

\begin{enumerate}

\item[1.]
Adding a pendant vertex to a vertex $v$, i.e. adding a new
vertex $v'$ and an edge $vv'$.

\item[2.]
Adding a twin vertex of a vertex $v$.

\begin{itemize}

\item[a.]
True twin:
adding a new vertex
$v'$ adjacent to $v$ and edge $uv'$ if and only if $uv$ is an edge.

\item[b.]
False twin: adding a new vertex $v'$, and edge $uv'$ if and only if
$uv$ is an edge, but not the edge $vv'$.

\end{itemize}

\end{enumerate}
\end{defi}

Note that a DH graph is a circle graph with operation $1$ corresponding
to adding a small chord perpendicular to $v$ in the chord diagram,
operation $2a$
corresponding to adding a very close parallel chord,
and $2b$ to adding a very close crossing chord as in Figure $5$.\\

[Insert Figure $5$: Effect of adding a pendant or twin
vertex.]\\

\begin{prop}(Bandelt and Mulder \cite{BM86}, Corollary $3$)

A DH graph is bipartite if and only if no true twins are added in
its construction.

\end{prop}

We abreviate bipartite distance hereditary graphs as BDH graphs.

The following proposition highlights a duality among pendant and
twin vertices with respect to the pivot operation.

\begin{prop}\label{5.7a}
If $uv$ is an edge of $G$, then $w$ is pendant on
$u$ in $G$ if and only if $w$ is a false twin of $v$ in $G^{uv}$,
and $u$, $v$ are true twins in $G$ if and only if
$u$, $v$ are true twins in $G^{uv}$.
\end{prop}
\begin{proof}
This follows because, if $w$ is pendant on $u$, then $w\in A_u$ but
is adjacent to no vertices in $A_{uv}$ or $A_v$
($A_u$, $A_w$, $A_{uw}$ as in Section \ref{complexity}). Thus,
in $G^{uv}$, $w$ has edges joining it to all of the neighbors of $v$, and
no others. Similarly, if $w$ is a false twin of $v$ in
$G^{uv}$, it loses all its edges except the one joining it
to $u$ in $(G^{uv})^{uv}=G$.
For true twins $u$, $v$, $A_u$ and $A_v$ are empty, so no toggling
occurs and in fact $G=G^{uv}$.
\end{proof}

\footnotetext{DH Graphs and the Interlace Polynomial.}

\begin{prop}\label{5.10}
If $G'$ is the graph that results from adding a
pendant vertex $w$ to a vertex $u$ of a loopless graph $G$, then
$q(G';x,y)=q(G;x,y)+ (x^2-2x+y) q(G- u;x,y)$.
\end{prop}
\begin{proof}
We pivot on edge $uw$ of $G'$, noting that
$G'- w=G$, and $G'^{uw}=G'$ since $w$ pendant on $u$ implies
that $A_w$ and $A_{uw}$ are empty, and
$G'- u=G- u \cup w$, where $w$ is an isolated
vertex. By (\ref{2eq}), since $G$, and hence $G'$, is loopless,
$q(G')=q(G'- w)+q(G'^{uw}- u)+
(x^2-2x)q(G'^{uw}- u- w)$
$=q(G)+q(G'- u)+ (x^2-2x)q(G'- u- w)$
$=q(G)+q(G- u\cup w)+ (x^2-2x)q(G- u)$
$=q(G)+(x^2-2x+y)q(G- u).$

\end{proof}

\begin{prop}\label{5.12}

If $G''$ is the graph that results from adding a false twin $w$ to
a vertex $v$ of a loopless graph $G$,
when $v$ is not isolated in $G$, then

$$q(G'';x,y)=q(G;x,y)+y(q(G;x,y)-q(G- v;x,y)).$$

\end{prop}
\begin{proof}

We pivot on edge $uv$ of $G''$, noting that since $v$  and $w$ are
false twins then $G''^{uv}$ is
$G^{uv}$ with vertex $w$ pendant on $u$. This follows because, except for $w$, $A_u$, $A_v$ and $A_{uv}$ are exactly
the same in both $G''$ and $G$, but $w$, adjacent precisely to everything in $A_v\cup A_{uv}$, loses
all edges except $wu$ in the pivot. Also note that since $v$ and $w$ are duplicates
$G''- v=G''- w=G$. Now by (\ref{2eq}), we have

$$q(G'')=q(G''- v)+q(G''^{uv}- u)+(x^2-2x)q(G''^{uv}- u- v)$$

$$=q(G)+q(G^{uv}- u \cup w)+(x^2-2x)q(G^{uv}- u-v\cup w)$$

$$=q(G)+yq(G^{uv}- u)+y(x^2-2x)q(G^{uv}- u- v)$$

$$=q(G)+y(q(G)-q(G- v)).$$

\end{proof}

\begin{prop}\label{M}

If $G'''$ is the graph that results from adding a true twin $w$ to
$v$ of a loopless graph $G$, when $v$ is not isolated in $G$, then

$$q(G''';x,y)=2q(G;x,y)+((x-1)^{2}-1)q(G-v;x,y).$$

\end{prop}
\begin{proof}

Recall that if $v$ and $w$ are true twins, then $G'''^{vw}=G'''$ and
hence $G'''^{vw}-w=G'''-w=G$, so the result follows from
equation (\ref{2eq}).

\end{proof}

The duality among pendant and twin vertices is particularly
apparent in the context of the vertex-nullity interlace polynomial:

\footnotetext{DH Graphs and the Interlace Polynomial.}

\begin{cor}\label{4.14}

The vertex-nullity interlace polynomial has the following duality identities:

\begin{enumerate}

\item
If $G'$ is the graph that results from adding a pendant vertex $w$ to a
vertex $u$ of $G$, then

$$q_N(G';x)=q_N(G;x)+xq_N(G-u;x)$$.

\item
If $G''$ is the graph that results from adding a false twin
$w$ to a non-isolated vertex
$v$ with $u$ adjacent to $v$, then

$$q_N(G'';x)=q_N(G;x)+xq_N(G^{uv}-u;x)$$.

\item
If $G'''$ is the graph that results from adding a true twin $w$ to a
non-isolated vertex $v$, then

$$q_N(G''';x)=2q_N(G;x).$$

\end{enumerate}

\end{cor}
\begin{proof}
Parts $1$, $2$, $3$ follow from equation (\ref{stateqN}) and
Propositions \ref{5.10}, \ref{5.12} and \ref{M}, respectively.
Part 2) is also a restating of Arratia, Bollob\'as and Sorkin
\cite{ABS} (Proposition $40$), using an expression just before the
final form given there.

\end{proof}

\begin{cor}\label{5.11}

$\gamma$ is invariant under pendant or false twins, and doubles for
true twins, as follows:

\begin{enumerate}

\item
If $G'$ is the graph that results from adding a pendant vertex to $G$ with
$|V(G)|\geq2$, then $\gamma (G')=\gamma (G)$, i.e. $\gamma$ is invariant
under the addition of pendant vertices.

\item
If $G''$ is the graph that results from adding a false twin to a
non-isolated
vertex of $G$, then $\gamma (G'')=\gamma (G)$, so
$\gamma$ is invariant under the addition of false twin vertices.

\item
If $G'''$ is the graph that results from adding a true twin to a
non-isolated vertex of $G$, then $\gamma (G''')=2\gamma (G)$, so
$\gamma$ is doubled by the addition of true twin vertices.

\end{enumerate}
\end{cor}

\footnotetext{DH Graphs and the Interlace Polynomial.}

\begin{defi}\label{5.8}
\rm
Let $G$ and $F$ be graphs with $u\in V(G)$, $v\in V(F)$. Then the
one point join of $G$ and $F$, denoted by
$G  _{u}\!\!\cdot _{v} F$ is
formed by identifying $u$ and $v$, resulting in a
cut vertex of $G _{u}\!\!\cdot _{v}  F$. Equivalently, if
$H$ is a graph with a cut vertex $v$, then $H$ is the
one point joint of $G$ and $F$ where
$G=v\cup G'$, for $G'$ a component of $H- v$ containing
at least one neighbor of $v$ in $H$, and
$F=v\cup F'$, for $F'$ the complement of $G'$ in
$H- v$.
\end{defi}

Note that adding a pendant edge to $G$ at $v$ is equivalent to taking
the one point join of $K_2$ to $G$ at $v$. Also note that if $v$
is an isolated vertex of $F$, then $G _u\cdot_{v} F$
is just the disjoint union of $G$ and $F- v$.

\footnotetext{DH Graphs and the Interlace Polynomial.}

\begin{prop}\label{5.14}
If $H$ is the one point join $G_u\cdot_v F$ where neither $u$ nor $v$ are
isolated vertices, then $2\gamma (H)=\gamma (G)\gamma (F)$.
\end{prop}
\begin{proof}
$H$ has more than one component if and only if at least one of $G$ or $F$
has, in which case both sides of the equation are zero. Thus, we may assume $H$ is
connected, and we proceed by induction on the number of vertices
of $G$. Since $u$ is not isolated and $G$ is connected, the base case is
$G=K_2$, i.e. a pendant vertex and the result follows from Corollary
\ref{5.11} and that $\gamma (K_2)=2$.

Now suppose $G$ has $n$ vertices. If every edge of $H$ is incident with
$u=v$, then $H=K_{1,r}$, $G=K_{1,n}$, $F=K_{1,s}$ with $r=n+s$.
Since by Arratia, Bollob\'as and Sorkin \cite{ABS00}, Proposition 7.1,
$\gamma (K_{1,m})=2$ for $m\geq 1$, the result follows.

Otherwise, without loss of generality, there is an edge $ab$ in $G$ with
$u\notin \{a,b\}$. Note that toggling in $H$
with respect to $ab$ occurs only among the
edges of $G$ not $F$, so $H- a=(G- a)_u\cdot _v F$
and   $H^{ab}- b=(G^{ab}- b)_u\cdot _v F$.
Thus, $\gamma (H)=\gamma ((G- a)_u\cdot _v F)+
\gamma ((G^{ab}- b)_u\cdot _v F)$
$=\frac{1}{2}\gamma (G)\gamma (F)$ by
induction. Therefore,
$2\gamma (H)=\gamma (G)\gamma (F)$.

\end{proof}

\begin{defi}\label{5.9}
\rm
Let $G$ and $F$ be graphs with  $u\in V(G)$ and
$v\in V(F)$. Then the two point join of $G$ and $F$, denoted
$G _{u}\!\!:_{v} F$ is formed by adding edge $au$ whenever $av\in E(F)$ and
edge $bv$ whenever $bu\in E(G)$.
\end{defi}

Note that the vertices $u$ and $v$ are false twins in
$G _{u}\!\!:_{v}  F$. Also
adding a false twin $u$ of $v$ to $G$ is equivalent to
taking the two point join $G_u\!\!:_v F$ where $F$ consists of just
the single vertex $v$.

\begin{prop}\label{5.15}
If $H$ is the two point join $G_u\!\!:_v F$, then
$2\gamma (H)=\gamma (G)\gamma (F)$.
\end{prop}
\begin{proof}
Since $u$ and $v$ are duplicate vertices in $H$, by Corollary \ref{5.11},
$\gamma (H)=\gamma (H- u)$ but
$H- u =G_u\cdot _v F$, so by Proposition
\ref{5.14}, $\gamma(H- u)=\frac{1}{2}\gamma (G)\gamma (F)$,
so $2\gamma (H)=\gamma (G)\gamma (F)$.
\end{proof}

For the following, we need to recall alternative characterizations of
BDH graphs, with parts i-iii due to Bandelt and Mulder
\cite{BM86}, and the equivalence of
iii and iv, where a $(6,2)$-chordal graph is a graph such that
every cycle of length at least $6$ has at least $2$
chords (see Ausiello, D'Atri, and  Moscarini \cite{AD'AM86}).


\begin{prop}\label{N}(Bandelt and Mulder \cite{BM86}, Corollaries $3$ and $4$)

The following are equivalent:

\begin{itemize}

\item[i.]
$G$ is a BDH graph,

\item[ii.]
$G$ is constructed from a single vertex by a sequence of adding
pendant vertices and false twins, but no true twins,

\item[iii.]
$G$ is triangle-free and does not contain $C_n$ for $n>4$, nor the graph
consisting of $C_6$ with a chord connecting two antipodal vertices,

\item[iv.]
$G$ is a bipartite $(6,2)$-chordal graph.

\end{itemize}

\end{prop}
\begin{proof}

The equivalence of cases i-iii appears in Bandelt and Mulder
\cite{BM86}, and that iii is
equivalent to iv comes from noting that bipartite implies triangle
free, and if $H$ is a cycle then it has at least two chords, so there
are no induced cycles nor a $C_6$ with and antipodal chord, and
thus iv implies iii. On the other hand, if $G$ satisfies iii, then it is
bipartite since iii is equivalent to i. If $C_n$ is a cycle in $G$ with
$n>6$, then it has at least one chord, which
creates a $C_m$ with $m>4$, so it has a chord, and thus $C_n$ has at least
two chords. For $C_6$, since there are no triangles, and it can't
have just the antipodal chord, it must have two chords, and thus
$G$ is bipartite $(6,2)$-chordal.

\end{proof}

\begin{cor}\label{5.17}

If $H=G  _{u}\!\!\cdot _{v} F$, then $H$ is a BDH graph if and only
if both $G$ and $F$ are BDH graphs.

\end{cor}
\begin{proof}

This follows immediately from the characterization of a BDH graph as
a bipartite $(6,2)$-chordal graph, since any cycle of $H$ must be
entirely contained in either $G$ or $F$.

\end{proof}

We are now ready to prove the main theorem of this section,
which completes the classification begun by Aigner and
van der Holst in \cite{AH04} of
graphs for which $\gamma =2$.

\footnotetext{DH  Graphs and the Interlace Polynomial.}

Note that in Theorem \ref{5.18} we require that $H$ be a
simple graph. In fact, if $H$ consists of $m\geq 2$
parallel edges, then $H$ is a connected  graph
with $\gamma (H)=2$ but $H$ is not a
BDH graph.

\begin{thm}\label{5.18}
$H$ is a simple, connected graph with
$\gamma (H)=2$ if and only if $H$ is a  BDH graph
with at least two vertices.
\end{thm}
\begin{proof}

If $H$ is a BDH graph, then by Proposition \ref{N}
it is simple, connected,
and constructed using only pendant
vertices and false twins. Thus, $\gamma (H)=2$ follows from $\gamma (K_2)=2$ and
Corollary \ref{5.11}.

If $H$ is a simple connected graph with $\gamma (H)=2$, we
proceed by induction on $|V(H)|$ to show that $H$ is a BDH
graph with at least two vertices.

\footnotetext{DH Graphs and the Interlace Polynomial.}

If $|V(H)|\leq 1$, then $\gamma (H)\neq 2$. Thus $H$ has at least
two vertices. Moreover $H$ is connected by Property \ref{5.3}. If
$|V(H)|=2$, and $\gamma (H)=2$ then, since $H$ is simple and
connected, $H=K_2$, a BDH graph. So now assume
$|V(H)|=n\geq 3$, and pivot on an edge $uv$ of $H$.
Now $2=\gamma (H)=\gamma (H- v)+\gamma (H^{uv}- u)$.
Since $|V(H)|\geq 3$, neither summand can be $1$ by Property
\ref{5.5}, so one must be $0$ and the other $2$
(neither can be negative by Proposition \ref{pos}). If
$\gamma (H- v)=0$, then, by Property \ref{5.4}, $v$ is a cut
vertex of $H$. Since $v$ is a cut vertex, $H-v$ has at least two
components, so there are subgraphs $H_1$ and $H_2$, each with
more than one vertex, such that
$H=(H_1)  _{v}\!\cdot _{v} (H_2)$, as in
Definition \ref{5.8}. Hence by property \ref{5.5}, neither
$\gamma (H_1)$ nor $\gamma (H_2)$ is equal
to $1$.
By Proposition \ref{5.14}, $2\gamma(H)=\gamma (H_1)\gamma (H_2)$,
so $\gamma (H_1)=\gamma (H_2)=2$ and hence by induction $H_1$ and
$H_2$ are BDH graphs, and by Corollary \ref{5.17} $H$ is a
BDH graph.

If $\gamma (H^{uv}- u)=0$, we similarly have that
$H^{uv}$ is a BDH graph. Let $w$ be the last vertex added in a
construction of $H^{uv}$. If none of $w$ or any of its
neighbors include $u$ or $v$, then $w$ is a pendant or
duplicate edge in $H$ as well, and hence by
Corollary \ref{5.11},
$2=\gamma (H)=\gamma (H- w)$, so by induction
$H- w$ is a BDH graph and hence $H$ is.

If $w$ is pendant on $u$ or $v$ in $H^{uv}$, by Proposition
\ref{5.7a}, $w$ duplicates $v$ or $u$ in $H$, so
$2=\gamma (H)=\gamma (H- w)$, and hence, $H$ is a
BDH graph as above.

Similarly, using Proposition \ref{5.7a}, if $w$ duplicates $u$ or
$v$ in $H^{uv}$, then $w$ is pendant on $v$ or $u$ in
$H$, and hence, $H$ is a BDH graph.

Since $uv$ is an edge, $u$ and $v$ cannot be duplicate
vertices. Thus, the only remaining case is that one
is  pendant on the other. Without loss of generality, say $u$ is
pendant on $v$. In this case $A_u$ and $A_{uv}$ are empty so
$H^{uv}=H$, and $H$ has a pendant vertex and hence
is a BDH graph by induction as above.
\end{proof}

For DH graphs that are not necessarily bipartite, we have the
following property.

\begin{thm}\label{4.21B}

If $G$ is a DH graph with at least two vertices, and $n$ true twins
are added in some construction sequence of $G$, then
$\gamma (G)=2^{n+1}$.

\end{thm}
\begin{proof}

This follows from $\gamma (K_2)=2$ and Corollary \ref{5.11}.

\end{proof}

\footnotetext{DH Graphs and the Interlace Polynomial.}

The converse does not hold however: $\gamma (G)$ a power of $2$
does not necessarily mean that $G$ must be a DH graph.
For example, $C_6$ is not a distance hereditary graph, having
no pendant vertices nor twins either true or false. However, from
Arratia, Bollob\'as, and Sorkin
\cite{ABS00}, $q_N(C_6)=4x+10x^2+2x^3$, so $\gamma (C_6)=4$, a power of
$2$.

\begin{cor}\label{C}

If $G$ is a DH graph, then all construction sequences for $G$ must have
the same number of vertices added as true twins.

\end{cor}

Although Corollary \ref{C} may also be shown readily by induction, we
include it as an example of how structural information may be
encoded by the interlace polynomial.

\section{Relation of  BDH graphs to
series-parallel graphs
and polynomial time computability}\label{PD}

We now characterize BDH graphs in terms of their relation to
series-parallel graphs. Since BDH graphs may be recognized in
polynomial time and the Tutte polynomial computed in
polynomial time for series-parallel graphs, this chatacterization,
together with the relation between the Tutte and vertex-nullity
polynomials, will allow us to conclude that the vertex-nullity
polynomial is polynomial time to compute for the class of
BDH graphs.

Given a $2$-face colored, $4$-regular planar graph, we call the graph constructed
by placing
a vertex in each black face and connecting vertices whose faces share a vertex in the
original graph the {\em black face graph}.
A {\em digon} is a graph consisting
of two vertices joined by two edges in parallel. A
series-parallel graph is constructed from a digon by repeatedly
adding  edges in parallel to an existing edge
or subdividing an existing edge.

\begin{prop}\label{5.20}

$H$ is a  BDH graph with at least two vertices if and only
if it is the circle graph of an Euler circuit in $\vec{G}_m$, where $G$ is a
series-parallel graph.

\end{prop}

\begin{proof}

We proceed by induction on the number $n$ of vertices of $H$, or equivalently,
the number of edges of $G$.   If $n = 2$, then $G$ is a digon, so
the result is immediate,
since both possible cycles in the oriented medial graph $\vec{G}_m$ give a chord diagram
with two intersecting chords corresponding to the BDH graph $K_2$.
Furthermore, the only $4$-regular Eulerian digraphs with Euler circuits
that give rise to such a chord diagram have digons as their black face graphs,
and thus, are oriented medial graphs for a series-parallel graph.

\footnotetext{DH Graphs and the Interlace Polynomial.}

Now suppose the proposition holds whenever there are $n - 1$ vertices in
$H$ and edges in $G$, and suppose $H$ is a BDH graph with $n$ vertices.
Let $v'$ be the last vertex added in some construction of $H$,
and let $v$ be the vertex $v'$ either twins or is pendant upon.
By induction, let $G'$ be a series-parallel graph such that $H- v'$
is the circle graph of some Euler circuit in $\vec{G}'_m$.  Adding $v'$ to
$H- v'$
adds a parallel or small perpendicular edge in the chord diagram, as in
Figure $5$.

The effect in $\vec{G}'_m$ is to insert a small digon at $v$, with its interior face colored
white or black, depending on whether the original Euler circuit followed the white or
black faces, respectively,  in the case $v$ is pendant,
or followed the black or white faces respectively
in the case $v$ is a false twin.
If the interior of the digon is black, the effect is adding an edge in
series to $G'$
 to get the desired series-parallel graph $G$, and if the digon is white,
the effect is adding an edge in parallel.  See Figure $6$.\\

[Insert Figure $6$: Configurations in the medial graph.]\\

Similarly, if $G$ is a series-parallel graph with $n$ edges,
the same construction in reverse yields the desired BDH graph.

\end{proof}

The motivation for the connection between BDH and
series-parallel graphs arises from the desire for a class of
graphs on which the vertex-nullity interlace polynomial would be
tractable. Theorem \ref{B} gives a relation between the
vertex-nullity interlace
and the Tutte polynomial via a medial graph construction,
and Oxley and Welsh show that
the Tutte polynomial is polynomial time to compute for
series-parallel graphs in \cite{dominicT}. Thus, we consider graphs that
arise as circle graphs of Euler circuits in the oriented medial
graphs of series-parallel graphs, seeking a characterization of
such class of graphs.

Recalling that $K_4$ is the excluded minor for series-parallel graphs
(see Duffin \cite{D65} and Oxley \cite{Oxl82}), observe that the
only Euler circuits of the oriented medial graph of
$K_4$ give rise to $C_6$ and $C_6$ with a single
antipodal chord as circle graphs. This suggests that these
two graphs should be prohibited, leading us to the
bipartite $(6,2)$-chordal graphs, and
hence BDH graphs.

We also note that BDH graphs are
$4$-closed. The notion of $k$-closure was introduced in
\cite{Sar98} and used in \cite{Sar} in the more general context of matroids.
Graphs
that are $4$-closed are characterized by their closed sets of rank up to
four (see \cite{Sar98}).

\footnotetext{DH Graphs and the Interlace Polynomial.}

\begin{prop}\label{P}

If $G$ is a BDH graph and $v$ is a pendant or false twin vertex, then
$G-v$ is also a BDH graph.

\end{prop}
\begin{proof}

By Proposition \ref{N}, $G$ is bipartite $(6,2)$-chordal graph, and
clearly if $v$ is a pendant vertex then $G-v$ is still a bipartite
$(6,2)$-chordal, and hence BDH, graph. If $v$ is a false twin of $u$ in
$G$, and $C_n$ with $n\geq 6$, is a cylce of $G-v$, then it is a
cycle in $G$, and hence has two chords. The graph $G-v$ clearly remains
bipartite, and is connected since $v$ as a twin vertex cannot be a
cut vertex. Thus $G-v$ is still a bipartite $(6,2)$-chordal,
and hence BDH, graph.

\end{proof}

Recall from Bandelt and Mulder \cite{BM86} (Corollary $1$)
that every DH graph $G$ with
at least $4$ vertices has at least two disjoint twin pairs, or a
twin pair and a pendant vertex, or at least two
pendant vertices. Also note that all connected graphs on $3$ or
fewer vertices are DH graphs.

We now give the following elementary greedy algorithm for
recognizing BDH graphs in polynomial time. There are certainly more
sophisticated and general recognition algorithms for
BDH graphs (see  work by
Cicerone and  Di Stefano in \cite{CS99a} and
\cite{CS99b} for example), but we use
the following simplistic approach in order to leverage
computability properties of series-parallel graphs.

\begin{cor}\label{PDpoly}

A BDH graph may be recognized, and a construction sequence
found, in polynomial time.

\end{cor}
\begin{proof}

We can identify pendant vertices in $O(n^2)$ steps by examining each vertex
to determine if it is adjacent to exactly one other vertex. (In the case
that the graph information is stored in an adjacency list rather than
matrix, this can be done in $O(n)$ steps). We can identify
duplicate vertices in $O(n^3)$ steps by comparing the neighbors of each
of the $O(n^2)$ pairs of vertices. We use this to successively find and
remove pendant or twin vertices. The original graph is a BDH
graph if and only if the graph that remains at the end of this
process is a single vertex. If the resulting graph is a single
vertex, then reversing the order of vertex removals gives a
BDH construction sequence for the original graph.

\end{proof}

\footnotetext{DH Graphs and the Interlace Polynomial.}

We now show that BDH graphs form a tractable class of graph for the
vertex-nullity interlace polynomial.

\begin{thm}\label{5.21}
Let $H$ be a BDH graph. Then the vertex-nullity interlace polynomial
$q_N(H;x)$ of $H$ can be calculated in polynomial time.
\end{thm}
\begin{proof}
Oxley and Welsh \cite{dominicT} have shown that
the Tutte polynomial $t(G;x,y)$ of a
series-parallel graph $G$ can be
calculated in polynomial time. In fact,
Noble \cite{noble} has shown that it can be calculated in a linear number of
multiplications involving $O(|V|)$ factors.
By Proposition \ref{5.20}, $H$ is the circle graph of an Eulerian
circuit of $\vec{G}_m$, where $G$ is a series-parallel
graph.
By Corollary \ref{PDpoly} we can find a construction
sequence for $H$ in polynomial time.
Construct $G$ using the construction sequence
of $H$ as in Proposition \ref{5.20}.
By Theorem \ref{B}, $q_N(H;x)=t(G;x,x)$. Therefore
$q_N(H;x)$ can be calculated in polynomial time.
\end{proof}

In closing we observe that there are a number of graph classes closely
related to DH graphs (see
Brandst\"adt, Le, and Spinrad \cite{BLS99} and McKee and McMorris
\cite{MM99}), many
with construction methods similar to those for DH graphs,
and further investigation of these classes in relation to the interlace
polynomial may well prove a fruitful area of research.

\begin{center}
{\bf\em Acknowledgements:}
\end{center}

We would like to thank an anonymous
referee for suggesting a number of productive
areas of investigation,
Dr. Greta Pangborn for
several helpful discussions, and
Dr. Terry McKee for introducing us to distance hereditary
graphs.

\end{document}